\documentclass[12pt]{article}
\usepackage[final]{epsfig}
\usepackage{graphics}
\usepackage{amsmath}
\usepackage{amsfonts}
\usepackage{latexsym}
\usepackage{amssymb}
\usepackage{graphicx}
\usepackage{url}
\usepackage{epstopdf}

\newtheorem{lemma}{Lemma}[section]

\newtheorem{remark}[lemma]{Remark}

\newtheorem{theorem}{Theorem}

\begin{document}
\newcommand{\eps}{{\varepsilon}}
\newcommand{\proofend}{$\Box$\bigskip}
\newcommand{\C}{{\mathbb C}}
\newcommand{\Q}{{\mathbb Q}}
\newcommand{\R}{{\mathbb R}}
\newcommand{\Z}{{\mathbb Z}}
\newcommand{\RP}{{\mathbf {RP}}}
\newcommand{\CP}{{\mathbf {CP}}}
\newcommand{\Tr}{\rm Tr}
\def\proof{\paragraph{Proof.}}

\title{Skewers}

\author{Serge Tabachnikov\footnote{
Department of Mathematics,
Penn State University,
University Park, PA 16802;
tabachni@math.psu.edu}
}

\date{}
\maketitle

\section{Introduction} \label{intro}

Two lines in 3-dimensional space are skew if they are not coplanar. Two skew lines share a common perpendicular line that we call their {\it skewer}. We denote the skewer of lines $a$ and $b$ by $S(a,b)$.\footnote{One can also define the skewer of two intersecting lines: it's the line through the intersection point,  perpendicular to both lines.}

\begin{figure}[hbtp]
\centering
\includegraphics[height=1.7in]{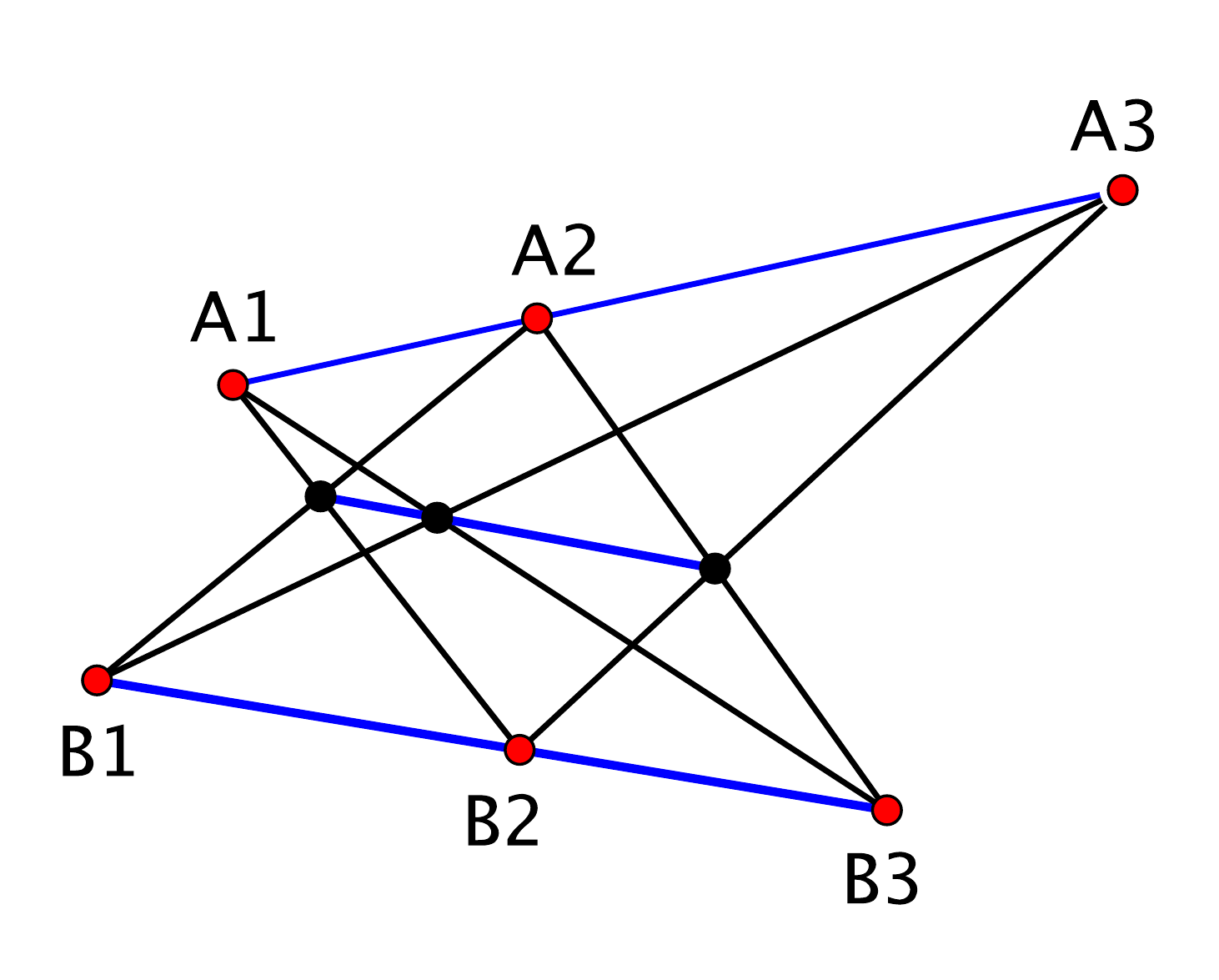}
\caption{The Pappus theorem.}
\label{Pappus}
\end{figure}

Consider your favorite configuration theorem of plane projective geometry that involves points and lines. For example, it may be the Pappus theorem, see Figure \ref{Pappus}: if $A_1,A_2,A_3$ and $B_1,B_2,B_3$ are two triples of collinear points, then the three intersection points $A_1B_2\cap A_2B_1$, $A_1B_3\cap A_3B_1$, and $A_2B_3\cap A_3B_2$ are also collinear (we refer to \cite{RG} for a modern viewpoint on projective geometry).

The Pappus theorem has a skewer analog in which both points and lines are replaced by lines in 3-space and the incidence between a line and a point translates as the intersection of  the two respective lines at  right angle. The basic 2-dimensional operations of connecting two points by a line or by intersecting two lines at a point translate as taking the skewer of two lines. 

\begin{theorem}[Skewer Pappus theorem I] \label{skPappus}
Let $a_1,a_2,a_3$ be a triple of lines with a common skewer, and let $b_1,b_2,b_3$ be another  triple of lines with a common skewer. Then the  lines 
$$
S(S(a_1,b_2),S(a_2,b_1)),\ S(S(a_1,b_3),S(a_3,b_1)),\  {\rm and}\ \ S(S(a_2,b_3),S(a_3,b_2))
$$
 share a skewer.
\end{theorem} 

In this  theorem,  we assume that the lines involved are in general position in the following sense: each time one needs to draw a skewer of two lines, this operation is well defined and unique. This assumption holds in a Zariski open subset of the set of the initial lines (in this case, two triples of lines with common skewers, $a_1,a_2,a_3$ and  $b_1,b_2,b_3$). A similar general position assumption applies to other theorems in this paper.\footnote{The configuration theorems of plane geometry also rely on similar general position assumptions.}

Another skewer analog of the Pappus theorem was discovered by R. Schwartz.

\begin{theorem}[Skewer Pappus theorem II] \label{othPapp}
Let $L$ and $M$ be a pair of skew lines. Choose a triple of points $A_1,A_2,A_3$ on $L$ and a triple of points $B_1,B_2,B_3$ on $M$. Then the lines 
$$
S((A_1 B_2), (A_2 B_1)), \  S((A_2 B_3), (A_3 B_2)),\  {\rm and}\ \ S((A_3 B_1), (A_1 B_3))
$$
share a skewer.
\end{theorem}

Although the formulation of Theorem \ref{othPapp} is similar to that of Theorem \ref{skPappus}, we failed to prove it along the lines of the proofs of other results in this paper, and the `brute force' proof of Theorem \ref{othPapp} is postponed until  Section \ref{Papprev}.

\begin{figure}[hbtp]
\centering
\includegraphics[height=1.7in]{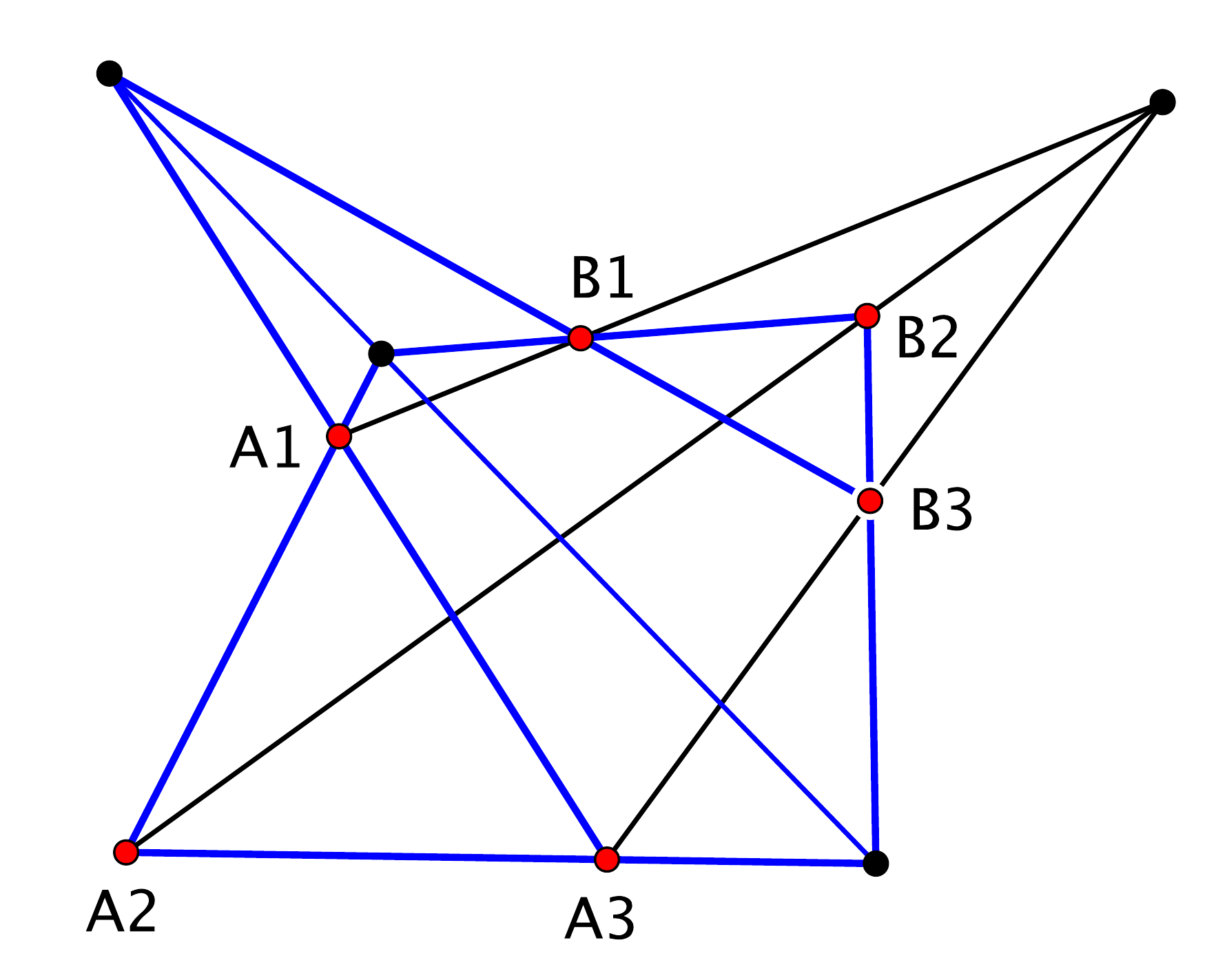}
\caption{The Desargues theorem.}
\label{Desargues}
\end{figure}

Another classical example is the Desargues theorem, see Figure \ref{Desargues}: if the three lines $A_1 B_1$, $A_2B_2$ and $A_3B_3$ are concurrent, then the three intersection points $A_1A_2\cap B_1B_2$, $A_1A_3\cap B_1B_3$, and $A_2A_3\cap B_2B_3$ are collinear. 

And one has a skewer version:

\begin{theorem}[Skewer Desargues theorem] \label{skDesargues}
Let $a_1,a_2,a_3$ and $b_1,b_2,b_3$ be two triples of lines such that the lines $S(a_1,b_1), S(a_2,b_2)$ and $S(a_3,b_3)$ share a skewer. Then the lines 
$$
S(S(a_1,a_2),S(b_1,b_2)),\ S(S(a_1,a_3),S(b_1,b_3)),\  {\rm and}\ \ S(S(a_2,a_3),S(b_2,b_3))
$$
also share a skewer.
\end{theorem}

The projective plane $\RP^2$ is the projectivization of 3-dimensional vector space $V$.
Assume that the projective plane is equipped with a polarity, a projective isomorphism $\varphi: \RP^2 \to (\RP^2)^*$ induced by a self-adjoint linear isomorphism $V \to V^*$. 

\begin{figure}[hbtp]
\centering
\includegraphics[height=1.5in]{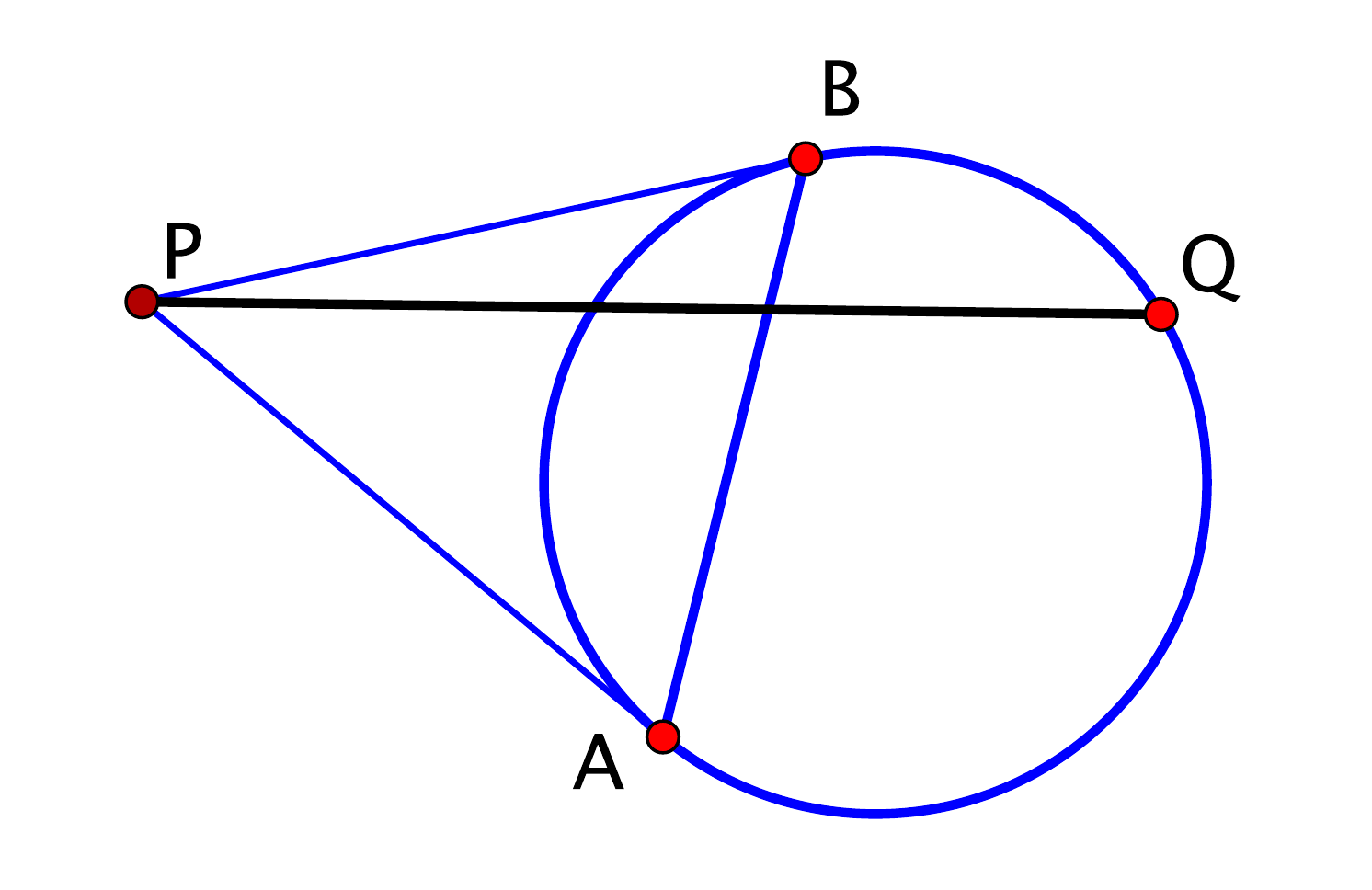}
\caption{Point $P$ is polar dual to the line $AB$.}
\label{polar}
\end{figure}

In particular, in 2-dimensional spherical geometry, polarity is the correspondence between great circles and their poles.\footnote{On $S^2$, this is a 1-1 correspondence between oriented great circles and points; in its quotient $\RP^2$, the elliptic plane, the orientation of lines becomes irrelevant.} In terms of  2-dimensional hyperbolic geometry, polarity is depicted in Figure \ref{polar}: in the projective model, $H^2$ is represented by the interior of a disc in $\RP^2$, and the polar points of lines lie outside of $H^2$, in the de Sitter world.

As a fourth example, consider a theorem that involves polarity, namely, the statement that the altitudes of a (generic) spherical or a hyperbolic triangle are concurrent (in the hyperbolic case, the intersection point may also lie in the de Sitter world). 

\begin{figure}[hbtp]
\centering
\includegraphics[height=2in]{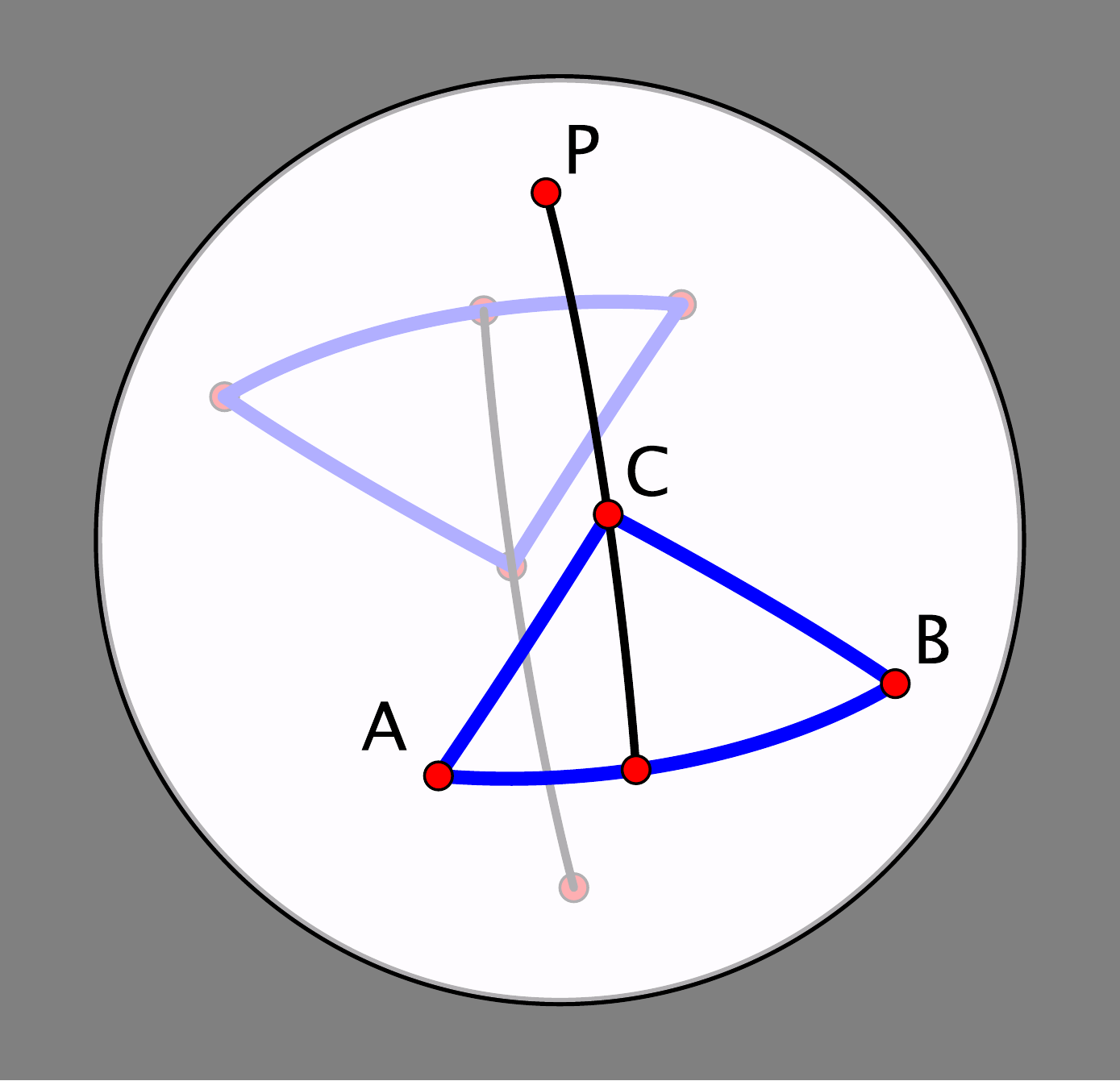}
\caption{Altitude of a spherical triangle.}
\label{sphere}
\end{figure}

The altitude of a spherical triangle $ABC$ dropped from vertex $C$ is the great circle through $C$ and the pole $P$ of the line $AB$, see Figure \ref{sphere}. Likewise, the  line $PQ$ in Figure \ref{polar} is orthogonal in $H^2$ to the line $AB$.  

In the skewer translation, we do not distinguish between polar dual objects, such as the line $AB$ and its pole $P$ in Figure \ref{sphere}. This yields the following theorem.

\begin{theorem}[Petersen-Morley  \cite{Mo}] \label{skAlt}
Given three lines $a,b,c$, the  lines
$$
S(S(a,b),c),\ S(S(b,c),a),\ \ {\rm and}\ \ S(S(c,a),b)
$$
share a skewer.\footnote{This result is also known as Hjelmslev-Morley theorem, see \cite{Fe}.}
\end{theorem}

\begin{figure}[hbtp]
\centering
\includegraphics[height=1.7in]{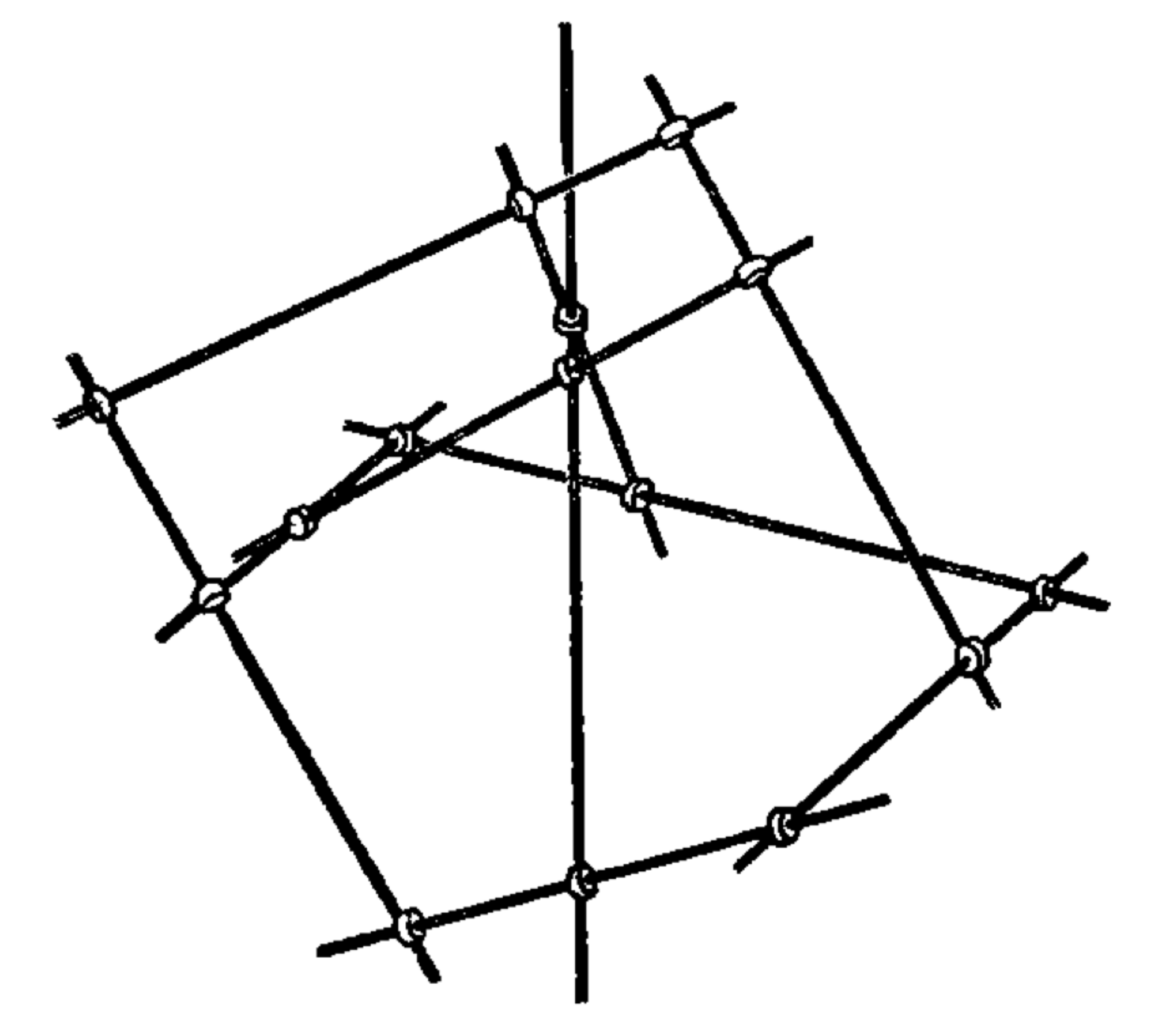}
\caption{Petersen-Morley configuration in Euclidean space.}
\label{ten}
\end{figure}

In words, {\it the common normals of the opposite sides of a rectangular hexagon have a common normal}; see Figure \ref{ten}, borrowed from \cite{Mo2}.

These `skewer' theorems hold not only in the Euclidean, but also in the elliptic and hyperbolic geometries. In $H^3$, two non-coplanar lines have a unique skewer. In  elliptic space $\RP^3$, a pair of  generic lines has two skewers; we shall address this subtlety in Section \ref{spherical}.

In the next section we shall formulate a general correspondence principle, Theorem \ref{principle}, establishing skewer versions of plane configuration theorems. This correspondence principle will imply  the above formulated theorems, except for Theorem \ref{othPapp}, whose proof will be given in Section \ref{Papprev}. 

The correspondence principle concerns line geometry of 3-dimensional projective space, a subject that was thoroughly studied in the 19th century by many an eminent mathematician (Cayley, Chasles, Klein, Kummer, Lie, 
Pl\"ucker, Study, ...) See \cite{Je} for a classical and \cite{PW} for a modern account.

Although we did not see the formulation of our Theorem \ref{principle} in the literature, we believe that  classical geometers would not be surprised by it. Similar ideas were expressed earlier. In the last section of \cite{Co}, 
H. S. M. Coxeter writes:
\begin{quote}
... every projective statement in which one conic plays a special role can be translated into a statement about hyperbolic space.
\end{quote}   
Coxeter illustrated this by the hyperbolic version of the Petersen-Morley theorem. 

Earlier F. Morley \cite{Mo1} also discussed the hyperbolic Petersen-Morley theorem, along with a version of Pascal's theorem for lines in $H^3$ (the ``celestial sphere" in the title of this  paper is the sphere at infinity of hyperbolic space).

We are witnessing a revival of projective geometry \cite{PW,RG}, not least because of  the advent of computer-based methods of study, including interactive geometry software (such as Cinderella\footnote{Which was used to create illustrations          
in this paper.} and GeoGebra). Elementary projective geometry has served as a source of interesting dynamical systems \cite{Sc1,Sc2}, and it continues to yield surprises \cite{ST}. We hope that this paper will contribute to the renewal of interest in this classical area. 

\bigskip
{\bf Acknowledgments}. 
The `godfather' of this paper is Richard  Schwartz whose question was the motivation for this project,  and who discovered Theorem  \ref{othPapp} and helped with its proof. I am grateful to Rich for numerous stimulating discussions on this and other topics. 

I am also grateful to I. Dolgachev, M. Skopenkov, V. Timorin, and O. Viro for their insights and contributions. Many thanks to A. Barvinok who introduced me to the chains of circles theorems.

I was supported by NSF grants DMS-1105442 and  DMS-1510055. Part of this work was done during my stay at ICERM; it is a pleasure to thank the Institute for the inspiring, creative, and friendly   atmosphere. 

\section{Correspondence principle} \label{2proofs}

\subsection{What is a configuration theorem?} \label{what}
We adopt the following `dynamic' view  of   configuration theorems. 

One starts with an initial data, a collection of labelled points $a_i$ and lines $b_j$ in $\RP^2$, such that, for some pairs of indices $(i,j)$, the point $a_i$ lies on the line $b_j$. 
One also has an ordered list of instructions consisting of two operations: draw a line through a certain pair of  points, or intersect a certain pair of  lines at a point. These new  lines and points also receive labels. 

The statement of a configuration theorem is that, among so constructed points and lines,  certain incidence relations hold, that is, certain points   lie on certain lines.   

Assume, in addition, that  a polarity $\varphi: \RP^2 \to (\RP^2)^*$ is given. We may think of lines in $\RP^2$ as points in $(\RP^2)^*$. The polarity takes one back to $\RP^2$, assigning the polar point to each line and vice versa. 

Given a polarity, one adds to the initial data that, for some pairs of indices $(k,l)$,  the point $a_k$ is polar dual to the line $b_l$.  One also adds to a list of instructions the operation of taking the polar dual object (point $\leftrightarrow$  line).
Accordingly, one adds to the statement of a configuration theorem that 
certain points are polar dual to certain lines. 

We assume that the conclusion of a configuration theorem holds for almost every initial configuration of points and lines satisfying the initial conditions, that is, holds for a Zariski open set of such initial configurations 
(this formulation agrees well with interactive  geometry software that makes it possible to perturb the initial data without changing its combinatorics).

In this sense, a configuration theorem is not the same as a configuration of points and lines as described in Chapter 3 of \cite{HC} or in \cite{Gr}: there, the focus is on whether a combinatorial incidence is realizable by 
points and lines in the projective plane. 

\begin{figure}[hbtp]
\centering
\includegraphics[height=2in]{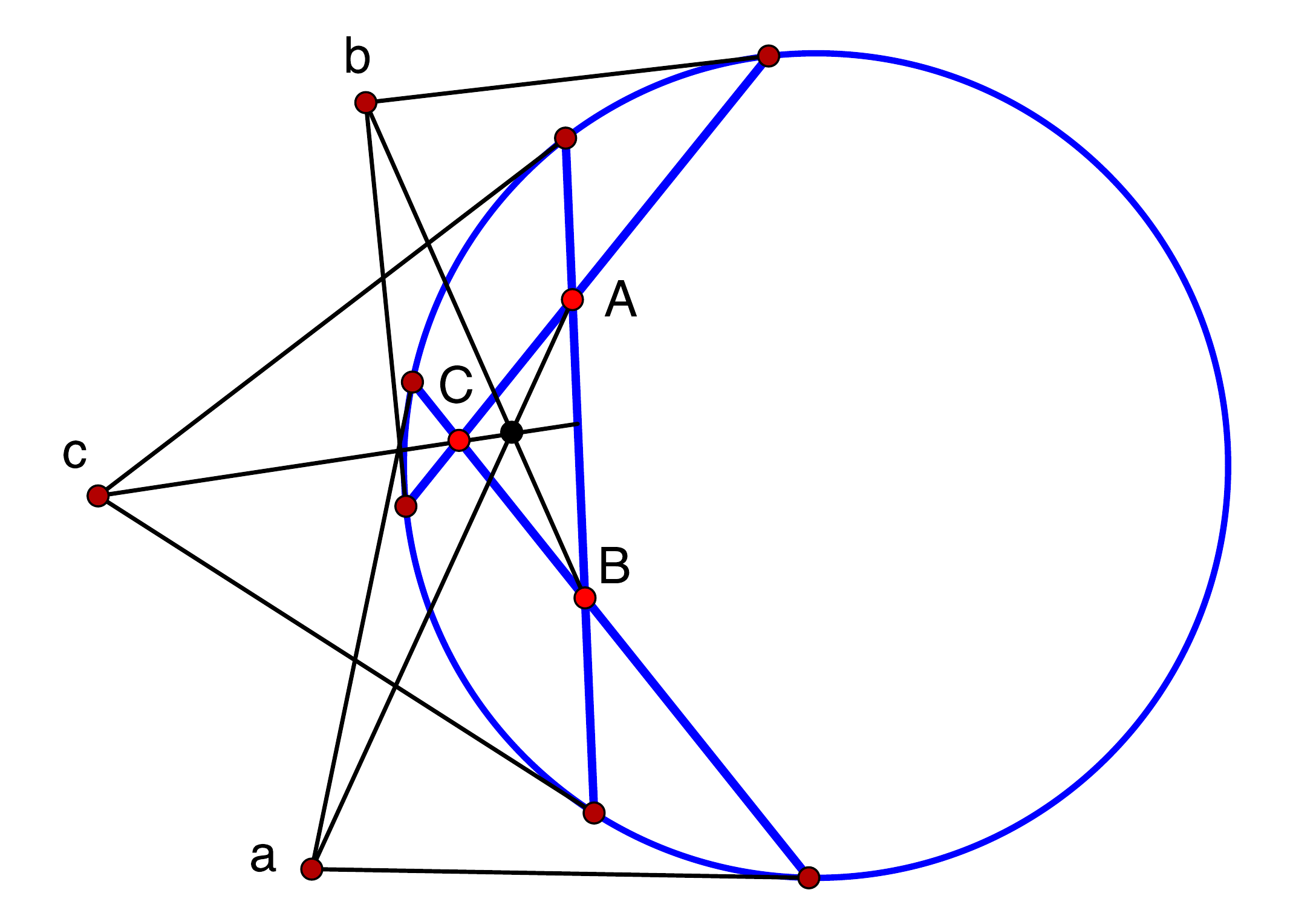}
\caption{Concurrence of the altitudes of a hyperbolic triangle.}
\label{hyptriangle}
\end{figure}

For example, the configuration theorem in Figure \ref{hyptriangle} has three points $A,B$ and $C$ as an initial data. One draws the lines $AB, BC$ and $CA$, and constructs their polar dual points $c, a$ and $b$, respectfully. Then one connects points $a$ and $A$, $b$ and $B$, and $c$ and $C$. The claim is that these three lines are concurrent (that is, the intersection point of the lines $aA$ and $bB$ lies on the line $cC$).

A configuration theorem for lines in space is understood similarly: one has an initial collection of labelled lines $\ell_i$ such that, for some pairs of indices $(i,j)$, the lines $\ell_i$ and $\ell_j$ intersect at right angle. There is only one operation, taking the skewer of  two lines. The statement of a configuration theorem is that certain pairs of thus constructed lines again intersect at right angle. This conclusion holds for almost all initial configurations of lines (i.e., a Zariski open set) satisfying the initial conditions. 

\subsection{Correspondence principle} \label{form}

The correspondence principle provides a dictionary that translates a plane configuration theorem, involving points and lines,   to a configuration theorem in space involving lines. 

\begin{theorem}[Correspondence principle] \label{principle}
To a plane configuration theorem with the initial data consisting of points $a_i$, lines $b_j$, and  incidences between them, there corresponds a configuration theorem for lines in space (elliptic, Euclidean, or hyperbolic), so that:  
\begin{itemize}
\item to each point $a_i$ and line $b_j$ of the initial data there corresponds a line in space;
\item whenever a point $a_i$ and a line $b_j$ are incident, the respective lines in space intersect at right angle;
\item the operations of  connecting two points by a line and of intersecting two lines at a point are replaced by the operation of taking the skewer of two lines.
\end{itemize}
If, in addition, a plane configuration theorem involves a polarity, then each pair of polar dual points and lines involved corresponds to the same line in space, and the operation of taking the polar dual object in the plane (point  $\leftrightarrow$ line) corresponds to the trivial operation of leaving a line in space intact. 
\end{theorem}

The reader might enjoy formulating the skewer version of the whole {\it hexagrammum mysticum}, the collection of results, ramifying the Pappus theorem, due to Steiner, Pl\"ucker, Kirkman, Cayley and Salmon; see \cite{CR1,CR2,Ho} for a modern treatment. 

We shall present two proofs of the Correspondence principle, one concerning the elliptic, and another the hyperbolic geometry. Either proof implies the Correspondence principle for the other two classical geometries: if a configuration theorem holds in the elliptic geometry, then it also holds in the hyperbolic geometry, and vice versa, by `analytic continuation'. And either non-zero curvature version implies the Euclidean one as a limiting case.

This analytic continuation principle is well known in geometry; we refer to  \cite{AP1,Pa} where it is discussed in detail.


\subsection{Elliptic  proof} \label{spherical}

A line in elliptic space $\RP^3$ is the projectivization of a 2-dimensional subspace of $\R^4$, and the geometry of lines in $\RP^3$ is the Euclidean geometry of 2-planes in $\R^4$.  The space of oriented lines is the  Grassmannian $G(2,4)$ of oriented 2-dimensional subspaces in $\R^4$. 

To every oriented line $\ell$ in $\RP^3$ there corresponds its dual oriented line $\ell^*$: the respective oriented planes in $\R^4$ are the orthogonal complements of each other (the orientation of the orthogonal complement is induced by the orientation of the plane and the ambient space). The dual lines are equidistant and they have infinitely many skewers. The preimage of a pair of dual lines in $S^3$ is a Hopf link.

The next lemma collects the properties of the Grassmannian $G(2,4)$ that we shall use. These properties are well known, see \cite{GW} for a detailed discussion.

\begin{lemma} \label{Grass1}
1) The Grassmannian is a product of two spheres:  $G(2,4)=S^2_-\times S^2_+$. This provides an identification of an oriented line in $\RP^3$ with a pair of points of the unit sphere $S^2$: $\ell\leftrightarrow (\ell_-,\ell_+)$. \\
2) The antipodal involutions of the spheres $S^2_-$ and $S^2_+$ generate the action of the Klein group $\Z_2\times\Z_2$ on the space of oriented lines. The action is generated by reversing the orientation of a line and by taking the dual line.\\
3) Two lines  $\ell$ and $m$ intersect at  right angle if and only if $d(\ell_-,m_-)=d(\ell_+,m_+)=\pi/2$, where $d$ denotes the spherical distance in $S^2$.\\
4) The set of lines that intersect $\ell$ at  right angle coincides with the set of lines that intersect $\ell$ and $\ell^*$.\\
5) A line $n$ is a skewer of lines $\ell$ and $m$ if and only if $n_-$ is a pole of the great circle $\ell_- m_-$, and $n_+$ is a pole of the great circle $\ell_+ m_+$.\\
6) A pair of generic lines has exactly two skewers (four, if orientation is taken into account), and they are dual to each other.
\end{lemma}

\proof Given two planes in $\R^4$, there are two angles, say $0\leq\alpha\leq\beta\le\pi/2$, between them: $\alpha$ is the smallest angle made by a line in the first plane with the second plane, and $\beta$ is the largest such angle. 

Recall the classical construction of Klein quadric (see, e.g., \cite{Do,PW}). Given an oriented plane $P$ in $\R^4$, choose a positive basis $u,v$ in $P$, and let $\omega_P$ be the bivector $u\wedge v$, normalized to be unit. In this way we assign to every oriented plane a unit decomposable element in $\Lambda^2 \R^4$. The decomposability condition $\omega\wedge\omega=0$ defines a quadratic cone in 
$\Lambda^2 \R^4$, and the image of the Grassmannian is the spherization of this cone (the Klein quadric is its projectivization).

Consider the star operator in $\Lambda^2 \R^4$, and let $E_-$ and $E_+$ be its eigenspaces with eigenvalues $\pm1$. These spaces are 3-dimensional, and  $\Lambda^2 \R^4=E_-\oplus E_+$. Let $S^2_{\pm}$ be the spheres of radii $1/\sqrt{2}$ in $E_{\pm}$. Then the bivector $\omega_P$ has the components in $E_{\pm}$ of lengths $1/\sqrt{2}$, and hence $G(2,4)=S^2_-\times S^2_+$. We rescale the radii of the spheres to unit. Thus an oriented plane $P$ becomes a pair of points $P_\pm$ of a unit sphere. 

Let us prove claim 2). Orientation reversing of a plane $P$ changes the sign of the bivector $\omega_P$ corresponding to the antipodal involutions of both spheres. Let $e_1,\ldots,e_4$ be an orthonormal basis in $\R^4$. Then the following vectors form bases of the spaces $E_{\pm}$:
$$
u_{\pm}=\frac{e_1\wedge e_2 \pm e_3\wedge e_4}{2},\ v_{\pm}=\frac{e_1\wedge e_3 \mp e_2\wedge e_4}{2},\ w_{\pm}=\frac{e_1\wedge e_4 \pm e_2\wedge e_3}{2}.
$$
Without loss of generality, assume that a plane $P$ is spanned by $e_1$ and $e_2$. Then $P^{\perp}$ is spanned by $e_3$ and $e_4$. Since $e_1\wedge e_2=u_+ + u_-, e_3\wedge e_4=u_+ - u_-$, the antipodal involution of $S^2_-$ sends $P$ to $P^{\perp}$.

Given two planes $P$ and $Q$, one has two pairs of points on $S^2$: $(P_-,Q_-)$ and $(P_+,Q_+)$. Let $\alpha$ and $\beta$ be the two angles between $P$ and $Q$. Then
$$
d(P_-,Q_-)=\alpha+\beta,\quad d(P_+,Q_+)=\beta-\alpha,
$$
see \cite{GW}.

In particular, $P$ and $Q$ have a nonzero intersection when $\alpha=0$, that is, when $d(P_-,Q_-)=d(P_+,Q_+)$. Likewise, $P$ and $Q$ are orthogonal when $\beta=\pi/2$. It follows that the respective lines  intersect at  right angle when $d(P_-,Q_-)=d(P_+,Q_+)=\pi/2$. This proves 3) and implies 5).

In terms of bivectors, two lines  intersect  if and only if 
$\omega_P  \cdot * \omega_Q =0,$
and they intersect at  right angle if, in addition,
$\omega_P \cdot \omega_Q =0$. 
Here dot means the dot product in $\Lambda^2 \R^4$ induced by the Euclidean metric.
The duality $\ell \leftrightarrow \ell^*$ corresponds to the star operator on bivectors. This implies 4). 

Finally, given two lines, $\ell$ and $m$, consider the distance between a point of $\ell$ and a point of $m$. This distance attains a minimum, and the respective line is a skewer of $\ell$ and $m$. By the above discussion, the skewers of lines $\ell$ and $m$ are the lines that intersect the four lines $\ell, \ell^*, m$ and $m^*$. 
This set is invariant under  duality and, by an elementary application of Schubert calculus (see, e.g., \cite{Do}), generically consists of two lines. This proves  6). 
\proofend

Thus taking the skewer of a generic pair of lines is a 2-valued operation. However, by the above lemma, the choice of the skewer does not affect the statement of the respective configuration theorem. 

One can also avoid this indeterminacy by
 factorizing the Grassmannnian $G(2,4)$ by the Klein group, replacing it by the product of two elliptic planes $\RP^2_-\times\RP^2_+$. In this way, we ignore orientation of the lines and identify dual lines with each other. As a result, a generic pair of lines has a unique skewer.  

Now to the Correspondence principle.

Given a plane configuration theorem, we realize it in the elliptic geometry: the initial data consists of points $a_i$ and lines $b_j$ in $\RP^2$ with some incidences between them, and the polarity in $\RP^2$ is induced by the spherical duality (pole $\leftrightarrow$ equator).

Let us replace the lines by their polar points.  Thus the initial data is  a collection of points $\{a_i, b_j^*\}$ in the projective plane such $d(a_i, b_j^*)=\pi/2$ when the point $a_i$ is incident with the line $b_j$. 

Likewise, instead of connecting two points, say $p$ and $q$, by a line, we take the polar dual point to this line, that is, the cross-product $p\times q$ of  vectors in $\R^3$, considered up to a factor. In this way, our configuration theorem will involve only points, and its statement is that certain pairs of  points are at distance $\pi/2$.

Take another the initial collection, $\{\bar a_i, \bar b_j^*\}$, and consider the collection of pairs $\{(a_i,\bar a_i), (b_j^*,\bar b_j^*)\}$ in $\RP^2_-\times \RP^2_+$. According to Lemma \ref{Grass1}, one obtains a configuration of  lines $\{\ell_i, \ell_j\}$ in elliptic space such that if a point $a_i$ is incident with a line $b_j$  then the corresponding lines $\ell_i$ and $\ell_j$ intersect at right angle. This is the initial data for the  skewer configuration theorem. By varying the generic choices of  $\{a_i, b_j^*\}$ and $\{\bar a_i, \bar b_j^*\}$ satisfying the initial incidences, we obtain a dense open set of initial configurations of lines $\{\ell_i, \ell_j\}$.

Likewise, the operations that comprise the configuration theorem (connecting pairs of points by lines and intersecting pairs of lines) become the operation of taking the skewer of a pair of lines, and the conclusion of the theorem is that the respective pairs of lines intersect at right angle.

\subsection{Hyperbolic proof} \label{hyperbolic}

In a nutshell, a skewer configuration theorem in 3-dimensional hyperbolic space  is a complexification of a configuration theorem in the hyperbolic plane. We use ideas of F. Morley \cite{Mo2} and V. Arnold \cite{Ar}.

Consider the 3-dimensional space of real binary quadratic forms $ax^2+2bxy+cy^2$ in variables $x,y$, equipped with the discriminant quadratic form $\Delta=ac-b^2$ and the respective bilinear form. We view the Cayley-Klein model of the  hyperbolic plane as the projectivization of the set $\Delta >0$, the circle at infinity being given by $\Delta=0$. The projectivization of the set $\Delta <0$ is the 2-dimensional de Sitter world.

Thus points of $H^2$ are elliptic (sign-definite) binary quadratic forms, considered up to a  factor. To a line in $H^2$ there corresponds its polar point that lies in the de Sitter world, see Figure \ref{polar}. Hence  lines in $H^2$ are hyperbolic (sign-indefinite) binary quadratic forms, also considered up to a  factor. 

Consider the standard area form $dx\wedge dy$ in the $x,y$-plane. The space of smooth functions is a Lie algebra with respect to the Poisson bracket (the Jacobian), and the space of quadratic forms is its 3-dimensional subalgebra $sl(2,\R)$. The following observations are made in \cite{Ar}.

\begin{lemma} \label{bracket}
A point is incident to a line in $H^2$ if and only if the corresponding quadratic forms are orthogonal with respect to the bilinear form $\Delta$. 
Given two points of $H^2$, the Poisson bracket of the respective elliptic quadratic forms is a hyperbolic one, corresponding to the line through these points. Likewise, for two lines in $H^2$, the Poisson bracket of the respective hyperbolic quadratic forms is an elliptic one, corresponding to the intersection point of these lines.
\end{lemma}

A complexification of this lemma also holds: one replaces $\RP^2$ by $\CP^2$, viewed as the projectivization of the space of quadratic binary forms (and losing the distinction between sign-definite and sign-indefinite forms). 
The conic $\Delta =0$ defines a polarity in $\CP^2$.

Lemma \ref{bracket} makes it possible to reformulate a configuration theorem involving points and lines in $H^2$ as a statement about the Poisson algebra of quadratic forms. For example, the statement that the three altitudes of a hyperbolic triangle are concurrent, see Figure \ref{sphere} right,  becomes the statement that the commutators
$$
\{\{f,g\},h\},\ \ \{\{g,h\},f\},\ \ {\rm and}\ \ \{\{h,f\},g\} 
$$
are linearly dependent, which is an immediate consequence of the Jacobi identity 
$$
\{\{f,g\},h\}+ \{\{g,h\},f\},+  \{\{h,f\},g\} =0
$$
in the Poisson Lie algebra.

Likewise, the Pappus theorem follows from the Tomihisa's identity 
$$
\{f_1, \{\{f_2, f_3\}, \{f_4, f_5\}\}\} + \{f_3, \{\{f_2, f_5\}, \{f_4, f_1\}\}\} + \{f_5, \{\{f_2, f_1\}, \{f_4, f_3\}\}\} = 0
$$
that holds in $sl(2,\R)$, see \cite{To}, and  also \cite{Ai,Iv,Sk} for this approach to configuration theorems.

Now consider 3-dimensional hyperbolic space $H^3$ in the upper halfspace model. The isometry group is $SL(2,\C)$, and the sphere at infinity is the Riemann sphere $\CP^1$. 

A line in $H^3$ intersects the sphere at infinity at two points, hence the space of (non-oriented) lines is the configuration space of unordered pairs of points, that is, the symmetric square of $\CP^1$ with the deleted diagonal. Note that  $S^2(\CP^1)=\CP^2$ (this is a particular case of the Fundamental Theorem of Algebra, one of whose formulations is that $n$th symmetric power of $\CP^1$ is $\CP^n$). Namely, to two points of the projective line one assigns the binary quadratic form having zeros at these points:
$$
(a_1:b_1,a_2:b_2) \longmapsto (a_1y-b_1x)(a_2y-b_2x).
$$
Thus a line in $H^3$ can be though of as a complex binary quadratic form up to a factor. 

 The next result is contained in \S 52 of \cite{Mo2}.

\begin{lemma} \label{Jacobian}
Two lines in $H^3$ intersect at  right angle if and only if the respective binary quadratic forms 
$f_i=a_i x^2 + 2 b_i xy + c_i y^2,\ i=1,2$, are orthogonal with respect to $\Delta$:
\begin{equation} \label{ort}
a_1c_2-2b_1b_2+a_2c_1=0.
\end{equation}
If two lines correspond to binary quadratic forms $f_i=a_i x^2 + 2 b_i xy + c_i y^2,\ i=1,2$, then their
skewer  corresponds to the Poisson bracket (the Jacobian)
$$
\{f_1,f_2\} = (a_1b_2-a_2b_1)x^2 + (a_1c_2-a_2c_1) xy + (b_1c_2-b_2c_1) y^2.
$$
\end{lemma} 

If $(a_1:b_1:c_1)$ and $(a_2:b_2:c_2)$ are homogeneous coordinates in the projective plane and the dual projective plane, then (\ref{ort}) describes the incidence relation between points and lines. In particular, the set of lines in $H^3$ that meet a fixed line at  right angle corresponds to a line in $\CP^2$.

Suppose a configuration theorem involving polarity is given in $\RP^2$. The projective plane with a conic provide the projective model of the hyperbolic plane, see Figures \ref{polar} and \ref{hyptriangle}, so the configuration in realized in $H^2$. Consider the complexification, the respective configuration theorem in $\CP^2$ with the polarity induced by $\Delta$. According to Lemma \ref{Jacobian}, this yields 
a configuration of lines in $H^3$ such that the pairs of incident points and lines correspond to pairs of lines intersecting at right angle. 
 
Another way of saying this is by way of comparing Lemmas \ref{bracket} and \ref{Jacobian}:
the relations in the Lie algebras $sl(2,\R)$ and $sl(2,\C)$ are the same, hence to a configuration theorem in $H^2$ there corresponds a skewer configuration theorem in $H^3$. 

\subsection{Euclidean picture} \label{Eucpict}

The following description of the Euclidean case is due to I. Dolgachev (private communication).

Add the plane at infinity to $\R^3$; call this plane $H$. A point of $H$ represents a family of parallel lines in $\R^3$. For a line $L$ in $\R^3$, let $q(L)=L\cap H$ be its direction, that is, the respective point at infinity.

One has a polarity in $H$ defined as follows. Let $A$ be a point in $H$. This point corresponds to a direction in $\R^3$. The set of orthogonal directions constitutes a line $A^*$ in $H$; this is the line polar to $A$.

\begin{lemma} \label{lineinf}
Let $L$ and $M$ be skew lines in $\R^3$. Then 
$$
q(S(L,M))= q(L)^* \cap q(M)^*.
$$
\end{lemma}

\proof
The direction $q(L)^* \cap q(M)^*$ is orthogonal to $L$ and to $M$, and so is the skewer $S(L,M)$. This implies the result. 
\proofend

Thus the skewer $S(L,M)$ is constructed as follows: find points $q(L)$ and $q(M)$ of the plane at infinity $H$, intersect their polar lines, and construct the line through   point $q(L)^* \cap q(M)^*$ that intersect $L$ and $M$. This line exists and is, generically, unique: it is the intersection of the planes through point $q(L)^* \cap q(M)^*$ and line $L$, and through point $q(L)^* \cap q(M)^*$ and line $M$. 

To summarize, a skewer configuration in $\R^3$ has a `shadow' in the plane $H$: to a line $L$ there corresponds the point $q(L)$ that is also identified with its polar line $q(L)^*$. In this way, the shadow of a skewer configuration is the respective projective configuration in the plane $H$. For example, both Theorems \ref{skPappus} and \ref{othPapp} become the usual  Pappus theorem in $H$.

\subsection{Odds and ends} \label{rmks}

1). {\it Legendrian lift}. One can associate a skewer configuration  in $\RP^3$ to a configuration in $S^2$ using contact geometry. 

A cooriented contact element in $S^2$ is a pair consisting of a point and a cooriented line through this point. The space of cooriented contact elements is $SO(3)=\RP^3$. We consider $\RP^3$ with its metric of constant positive curvature (elliptic space). The projection $\RP^3\to S^2$ that sends a contact element to its foot point is a Hopf fibration. 

The space of contact elements carries a contact structure generated by two tangent vector fields: $u$ is the rotation of a contact element about its foot point, and $v$ is the motion of the foot point along the respective geodesic. The fields $u$ and $v$ are orthogonal to each other.

 A curve tangent to the contact structure is called Legendrian. A smooth cooriented curve in $S^2$ has a unique Legendrian lift: one assigns to a point of the curve the tangent line at this point.

Consider a configuration of points and (oriented) lines  in $S^2$. One can lift each point as a Legendrian line in $\RP^3$, consisting of the contact elements with this  foot point. Likewise, one can lift each line as a Legendrian line, consisting of the contact elements whose foot point lies on this line. As a result, a configuration of lines and points in $S^2$ lifts to a configuration of lines in $\RP^3$ intersecting at right angle, as described in Theorem \ref{principle}. 

The family of (oriented) Legendrian lines in $\RP^3$ is 3-dimensional; it forms the Lagrangian Grassmannian $\Lambda(2) \subset G(2,4)$.   In the classical terminology, the 3-parameter family of Legendrian lines in projective space is the null-system, \cite{Je,Do}.
\smallskip

2). {\it Comparing the elliptic and hyperbolic approaches}. The approaches of Sections \ref{spherical} and \ref{hyperbolic} are parallel. The sphere $S^2$ in Section \ref{spherical} is the spherization of $\R^3=so(3)$, the Lie bracket being the cross-product of vectors. The pole of a line $uv$ in $S^2$ corresponds to the vector $u\times v$ in $\R^3$. Thus the operations of connecting two points by a line and  of intersecting two lines are encoded by the Lie bracket of $so(3)$.

Likewise, the Poisson bracket of two quadratic forms in Section \ref{hyperbolic} can be identified with the Minkowski cross-product that encodes the  operations of connecting two points by a line and  of intersecting two lines. 

Note that $so(3)$ is the Lie algebra of motions of $S^2$, whereas $sl(2,\R)$ is the Lie algebra of motions of $H^2$, and the complex forms of these Lie algebras coincide. Interestingly, this Lie algebraic approach to configuration theorems fails in the Euclidean plane, see \cite{Iv} for a discussion; however, Euclidean skewer configurations, such as the Petersen-Morley theorem, can be described in terms of the Lie algebra of motions of $\R^3$, see \cite{Sk}. 

In both proofs, one goes from the Lie algebra of motions in dimension 2 to that in dimension 3. In the elliptic situation, we have $so(4)=so(3)\oplus so(3)$, and in the hyperbolic situation, the Lie algebra of motions of $H^3$ is $sl(2,\C)$. Accordingly, an elliptic skewer configuration splits into the product of two configurations in $S^2$, and a hyperbolic skewer configuration is obtained from a configuration in $H^2$ by complexification.
\smallskip

3). {\it Skewers in $\R^3$ via dual numbers}. One can approach skewer configurations in $\R^3$ using Study's dual numbers \cite{St}; see \cite{PW} for a modern account. 

Dual numbers are defined similarly to complex numbers: 
$$
a+\varepsilon b,\ {\rm where}\ a,b\in\R,\ {\rm and}\ \varepsilon^2=0.
$$
Dual vectors are defined analogously.

To an oriented line $\ell$ in $\R^3$ one assigns the dual vector $\xi_{\ell}=u+\varepsilon v$, where $u\in S^2$ is the unit directing vector of $\ell$, and $v$ is the moment vector: $v=P\times u$ where $P$ is any point of $\ell$. The vectors $\xi_{\ell}$ form the Study sphere: $\xi_{\ell}\cdot \xi_{\ell}=1$.

This construction provides an isomorphism between the isometry group of $\R^3$ and the group of dual spherical motions. Two lines $\ell$ and $m$ intersect at  right angle if and only if $\xi_{\ell}\cdot \xi_{m}=0$.
Thus skewer configurations in $\R^3$ correspond to configurations of lines and points in the Study sphere whose real part are the respective configurations in $S^2$.

\section{Circles}

Denote the set of  lines in 3-space that share a skewer $\ell$ by ${\cal N}_\ell$.
We saw in Section \ref{2proofs} that ${\cal N}_\ell$ is  an analog of a line in the plane.
Two-parameter families of lines in 3-space are called congruences. ${\cal N}_\ell$ is a linear congruence: it is the intersection of the Klein quadric with a 3-dimensional subspace $\RP^3 \subset \RP^5$, that is, it is defined by two linear equations in Pl\"ucker coordinates.

Now we describe line analogs of circles.

Let $\ell$  be an oriented  line in 3-space (elliptic, Euclidean, or hyperbolic). Let $G_\ell$ be the subgroup of the group of orientation preserving isometries that preserve $\ell$. This group is 2-dimensional. Following \cite{Ri}, we call the orbit $G_\ell(m)$ of an oriented line $m$ an {\it axial congruence} with $\ell$ as axis.

In particular, ${\cal N}_\ell$ is an axial congruence.

In $\R^3$ (the case considered in \cite {Ri}), the lines of an axial congruence with axis $\ell$ are at equal distances $d$ from $\ell$ and make equal angles $\varphi$ with it. One defines the dual angle between two oriented lines  $\varphi + \varepsilon d$, see \cite{PW}. The dual angle between the lines of an axial congruence and its axis is constant. 

Thus, in $\R^3$, an axial congruence consists of a regulus (one family of ruling of a hyperboloid of one sheet) and its parallel translations along its axis.

Likewise, one defines a complex distance between oriented lines $\ell$ and $m$ in $H^3$. Let $d$ be the  distance from $\ell$ to $m$ along their skewer $S(\ell,m)$, and let $\varphi$ be the angle between $m$ and the line $\ell'$, orthogonal to $S(\ell,m)$ in the plane spanned by $\ell$ and $S(\ell,m)$, and intersecting $m$. (Both $d$ and $\varphi$  have signs determined by a choice of orientation of the skewer). Then the complex distance is given by the formula $\chi(\ell,m)=d+i\varphi$, see \cite{Ma}. Again, the complex distance between the lines of an axial congruence and its axis is constant.

If $\ell_{1,2}$ and $m_{1,2}$ are the respective points on the sphere at infinity $\CP^1$ then 
$$
\cosh^2\left(\frac{\chi(\ell,m)}{2}\right)=[\ell_1,m_1,m_2,\ell_2],
$$
where the cross-ratio is given by the formula
$$
[a,b,c,d]=\frac{(a-c)(b-d)}{(a-d)(b-c)},
$$
see \cite{Ma}.

In the next lemma, $\CP^1$ is the `celestial sphere', that is, the sphere at infinity of $H^3$. 

\begin{lemma} \label{round}
Let $\psi:\CP^1\to\CP^1$ be a M\"obius (projective) transformation having two distinct fixed points. The family of lines connecting point $z\in \CP^1$ with the point $\psi(z)$ is an axial congruence, and all axial congruences are obtained in this way.
\end{lemma}

\proof Without loss of generality, assume that the fixed points of $\psi$ are $0$ and $\infty$, and let $\ell$ be the line through these points. Then $\psi(z)=cz$ for some constant $c\in\C$. One has
$
[0,z,cz,\infty]=[0,1,c,\infty]=c/(c-1).
$
Hence, for the lines $m$ connecting $z$ and $\psi(z)$, the complex distance $\chi(\ell,m)$ is the same.

Conversely, given an axial congruence, we may assume, without loss of generality, that its axis $\ell$ connects $0$ and $\infty$. Then $G_{\ell}$ consists of the transformations $z\mapsto kz,\ k\in\C$. Let $m$ be the line connecting points $w_1$ and $w_2$. Then the axial congruence $G_{\ell}(m)$ consists of the lines connecting points $k w_1$ and $kw_2=\psi(kw_1)$, with $\psi: z\mapsto (w_2/w_1) z$.
\proofend

In $S^3$, an axial congruence is characterized by the condition that the angles $\alpha$ and $\beta$ (see the proof of Lemma \ref{Grass1}) between the axis and the lines of the congruence are constant. It follows from the proof of Lemma \ref{Grass1} that an axial congruence is a torus, a product of  circles, one in $S^2_-$ and another in $S^2_+$.

Thus an axial congruence of lines is an analog of a circle in 2-dimensional geometry. The arguments from Section \ref{spherical} imply  analogs of the basic properties of circles:
\begin{enumerate}
\item If two generic axial congruences  share a line then they share a unique other line.
\item Three generic oriented lines belong to a unique axial congruence.
\end{enumerate}
(A direct proof of the first property: if the axes of the congruences are $\ell_1$ and $\ell_2$, and the shared line is $m$, then the second shared line is obtained from $m$ by reflecting in $S(\ell_1,\ell_2)$ and reverting the orientation).

Using the approach of Section \ref{2proofs}, one extends the Correspondence principle to theorems involving circles. For example, one has

\begin{figure}[hbtp]
\centering
\includegraphics[height=1.8in]{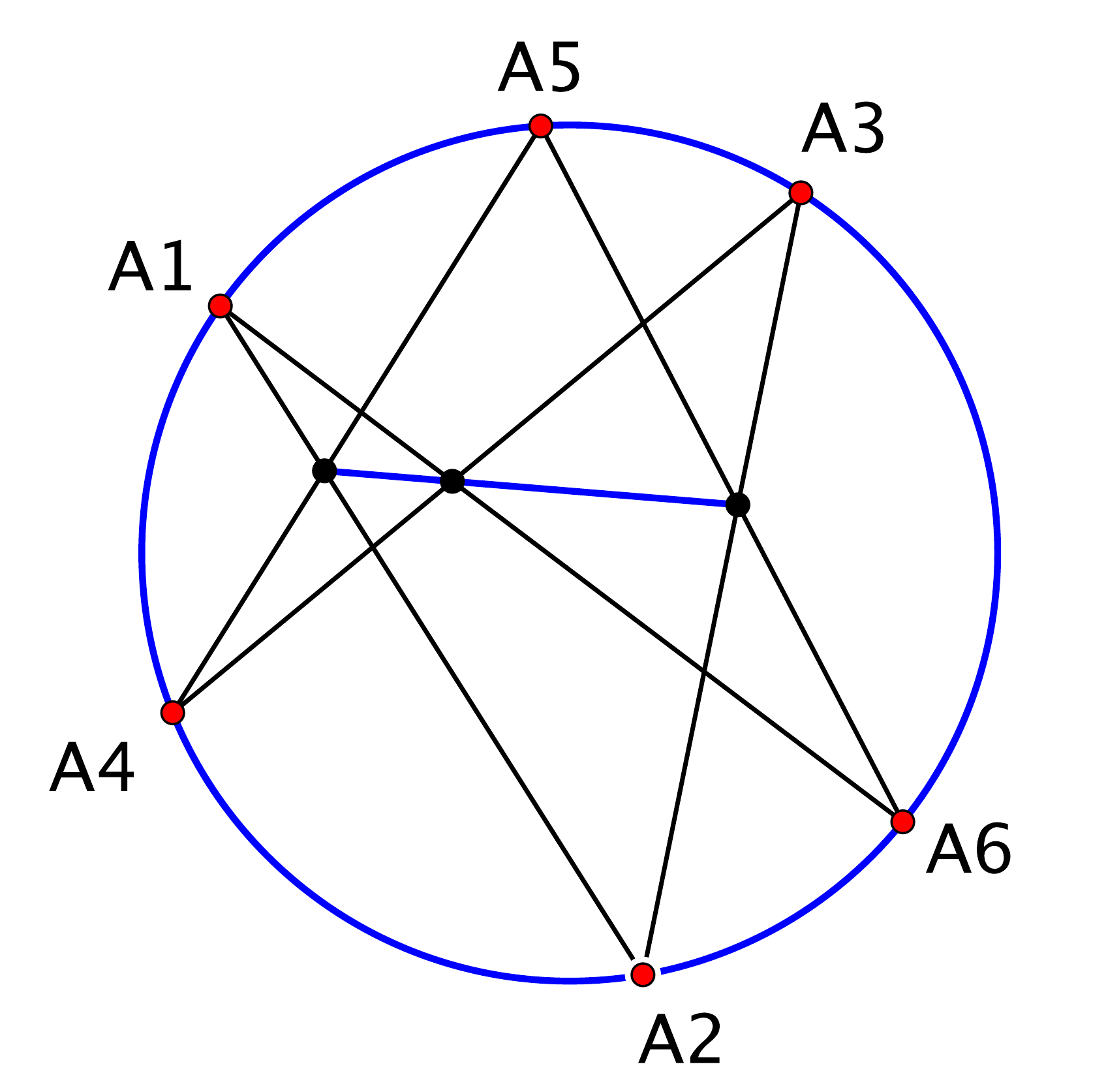}
\caption{Pascal's theorem for a circle.}
\label{Pascal}
\end{figure}

\begin{theorem}[Skewer Pascal theorem] \label{skPascal}
Let $A_1,\ldots,A_6$ be lines from an axial congruence. Then
$$
S(S(A_1,A_2),S(A_4,A_5)),\ S(S(A_2,A_3),S(A_5,A_6)),    \  {\rm and}\  S(S(A_3,A_4),S(A_6,A_1))
$$
share a skewer, see Figure \ref{Pascal}.  
\end{theorem}

\begin{figure}[hbtp]
\centering
\includegraphics[height=3in]{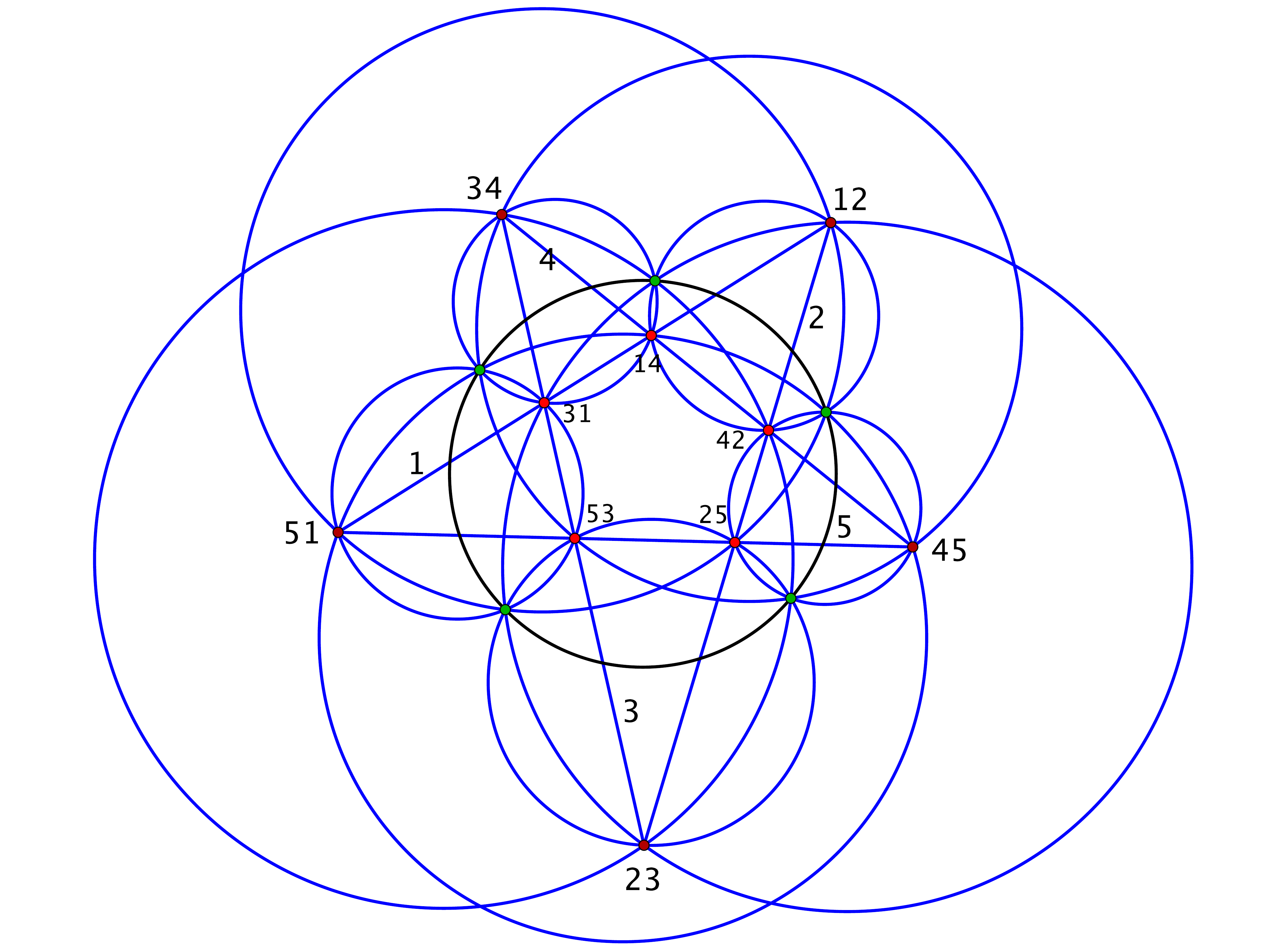}
\caption{Clifford's Chain of Circles ($n=5$).}
\label{Clifford}
\end{figure}

As another example, consider the  Clifford's Chain of Circles. This chain of theorems starts with a number of concurrent circles labelled $1,2,3,\ldots, n$. In Figure \ref{Clifford}, $n=5$, and the initial circles are represented by straight lines (so that their common point is at infinity).\footnote{As usual, lines are  considered as circles of infinite radius.} The intersection point of circles $i$ and $j$ is labelled $ij$. The circle through points $ij, jk$ and $ki$ is labelled $ijk$. 

The first statement of the theorem is that the circles $ijk, jkl, kli$ and $lij$ share a point;  this point is labelled $ijkl$. The next statement is that the points $ijkl, jklm, klmi, lmij$ and $mijk$ are cocyclic; this circle is labelled $ijklm$. And so on, with the claims of being concurrent and cocyclic alternating; see \cite{Coo,Mo2}, and \cite{KS,SK} for a relation with completely integrable systems.

A version of this theorem for lines in $\R^3$ is due to Richmond \cite{Ri}. The approach of Section \ref{2proofs} provides an extension to the elliptic and hyperbolic geometries.

\begin{theorem}[Clifford's Chain of Lines] \label{skClifford}
1) Consider axial congruences ${\cal C}_i,\ i=1,2,3,4$,  sharing a line.  For each pair of indices $i,j \in \{1,2,3,4\}$, denote by $\ell_{ij}$ the  line shared by ${\cal C}_i$ and ${\cal C}_j$, as described in statement 1 above. For each triple of indices $i,j,k \in \{1,2,3,4\}$, denote by ${\cal C}_{ijk}$  the axial congruence containing the lines $\ell_{ij},\ell_{jk},\ell_{ki}$, as described in the statement 2. Then the congruences ${\cal C}_{123}, {\cal C}_{234}, {\cal C}_{341}$ and ${\cal C}_{412}$ share a line. \\
2) Consider axial congruences ${\cal C}_i,\ i=1,2,3,4,5$,  sharing a line. Each four of the indices determine a line, as described in the previous statement of the theorem. One obtains five lines, and they all belong to an axial congruence.\\
3) Consider axial congruences ${\cal C}_i,\ i=1,2,3,4,5,6$,  sharing a line. Each five of them determine an axial congruence, as described in the previous statement of the theorem. One obtains six axial congruences, and they all share a line.
And so on...
\end{theorem} 

\begin{figure}[hbtp]
\centering
\includegraphics[height=1.7in]{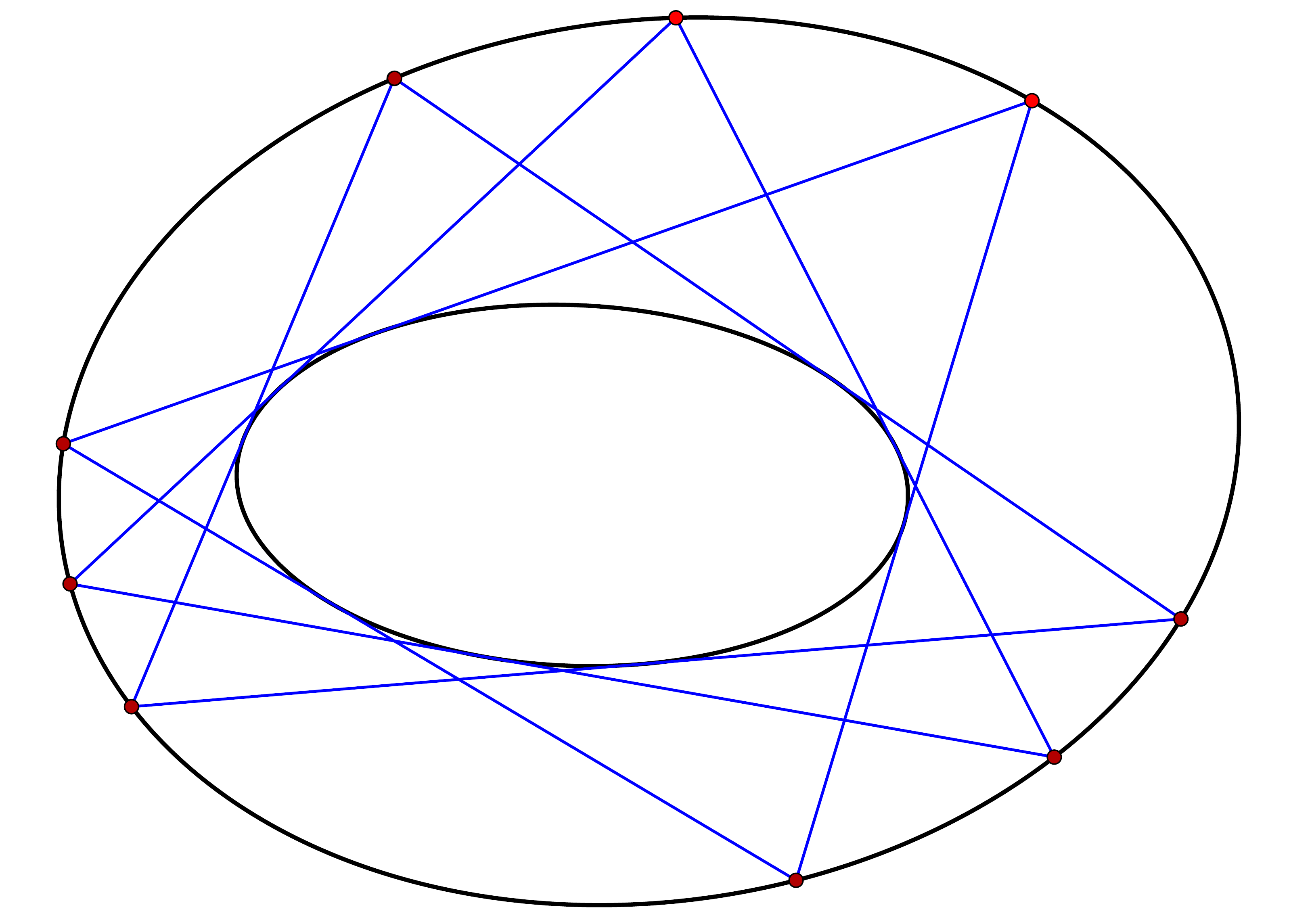}
\caption{Poncelet Porism, $n=3$.}
\label{Poncelet}
\end{figure}

Next, we present an analog of the Poncelet Porism, see, e.g., \cite{DR,Fl}. This theorem states that if there exists an $n$-gon inscribed into a conic and circumscribed about a nested conic  then every point of the outer conic is a vertex of such an $n$-gon, see Figure \ref{Poncelet}.

Consider a particular case when both conics are circles (a pair of nested conics can be sent to a pair of circles by a projective transformation).
The translation to the language of lines in space is as follows. 

Consider two generic axial congruences ${\cal C}_1$  and  ${\cal C}_2$, and assume that there exist a pair of lines $\ell_1\in {\cal C}_1$ and $\ell_2\in {\cal C}_2$ that intersect at  right angle. That is, ${\cal C}_1$ and ${\cal N}_{\ell_2}$ share the line $\ell_1$. 
By property 1) above, there exists a unique other  line $\ell_1'\in {\cal C}_1$, shared with ${\cal N}_{\ell_2}$, that is,  $\ell_1'$ intersects  $\ell_2$ at  right angle. Then there exists a unique other line $\ell_2'\in {\cal C}_2$ that intersects  $\ell_1'$ at  right angle, etc. We obtain a chain of intersecting orthogonal lines, alternating between the two axial congruences. 

The following theorem holds in the three classical geometries.

\begin{theorem}[Skewer Poncelet theorem] \label{skPoncelet}
 If this chain of lines closes up after $n$ steps, then the same holds for any starting pair of lines from ${\cal C}_1$  and  ${\cal C}_2$ that intersect at  right angle.
\end{theorem} 

\proof
Arguing as in Section \ref{spherical}, we interpret one axial congruence as the set of points of a spherical circle, and another one as the set of geodesic circles tangent to a spherical circle. The incidence between a geodesic and a point corresponds to two lines in space intersecting at  right angle. Thus the claim reduces to a version of the Poncelet theorem in $S^2$ where a spherical polygon is inscribed in a spherical circle and circumscribed about a spherical circle.

This spherical version of the Poncelet theorem is well known, see, e.g., \cite{CS,Ve}. For a proof, 
the central projection sends a pair of disjoint circles to a pair of nested conics in the plane, and the geodesic circles to straight lines, and the result follows from  the  plane Poncelet theorem.
\proofend

A pair of nested circles in the Euclidean plane is characterized by three numbers: their radii, $r<R$, and the distance between the centers, $d$. The conditions for the existence of an $n$-gons inscribed into one and circumscribed about another circle (a bicentric $n$-gon) are known as the Fuss relations. The first ones, for $n=3$ and $n=4$, are
$$
R^2-d^2=2rR,\quad  (R^2-r^2)^2=2r^2 (R^2+d^2);
$$
the case $n=3$ is due to Euler; Fuss found the relations for $n=4,\ldots,8$. More generally, Cayley gave conditions for Poncelet polygons to close up after $n$ steps for a pair of conics, see \cite{DR,Fl}).

It would be interesting to find an analog of the Fuss and Cayley relations; up to isometry, a pair of axial congruence depends on 6 parameters (two characterizing each congruence and two describing the mutual position of the axes). 

\section{Projections and conics} \label{conics}
In this section, we propose a definition-construction of a skewer analog of a conic. 

\begin{figure}[hbtp]
\centering
\includegraphics[height=2in]{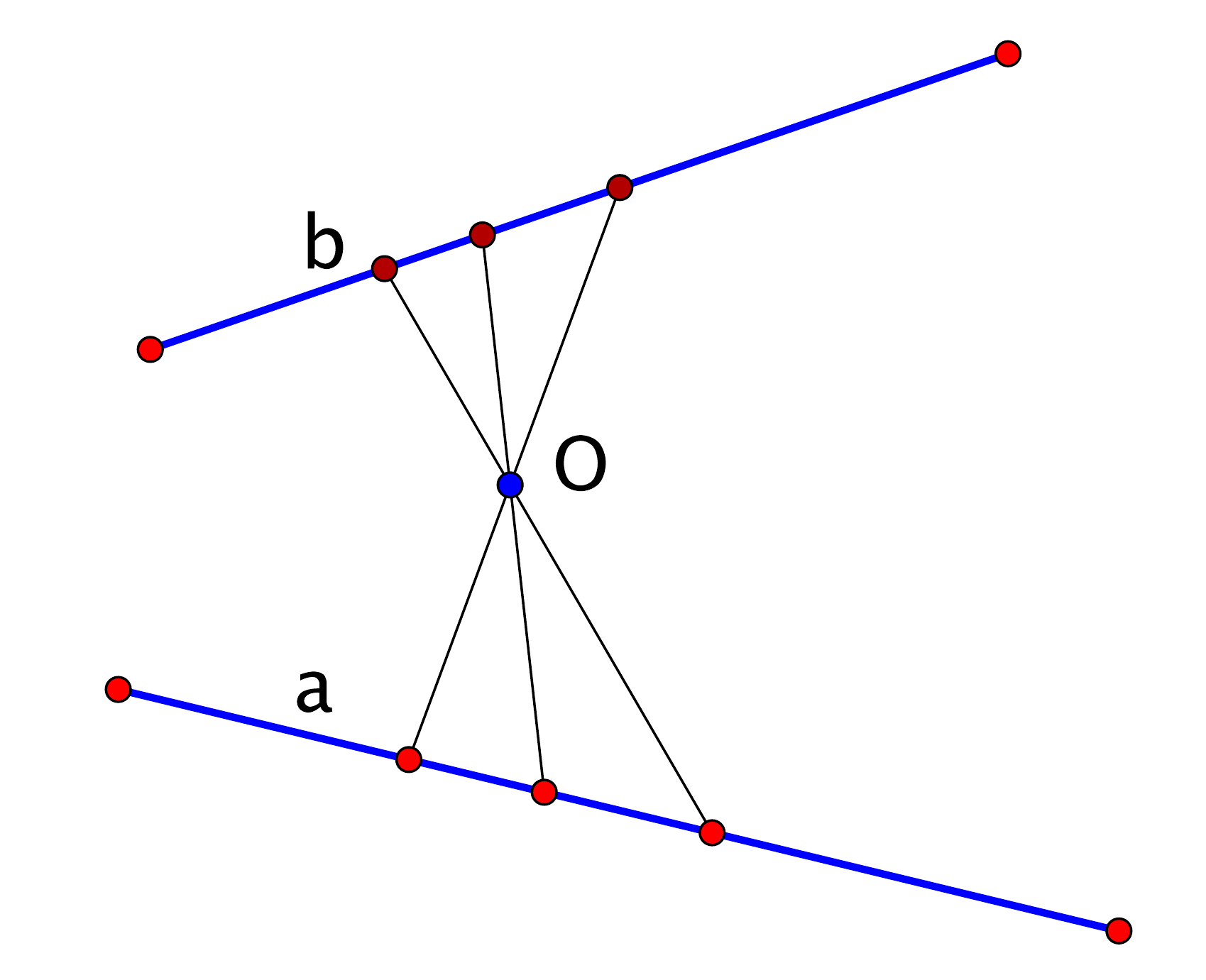}
\caption{A central projection of a line to a line.}
\label{projective}
\end{figure}

Let us first describe a skewer analog of a projection of a line to a line. Figure \ref{projective} depicts the central projection $\varphi_O: a \to b$ between two lines in the plane.

Consider three lines in space, $a,b$ and $O$, and define a map  $\varphi_O: {\cal N}_a \to {\cal N}_b$ as follows: for $\ell \in {\cal N}_a,$ set $\varphi_O (\ell) = S(S(\ell,O),b)$. This is a skewer analog of the central projection. Like in the plane, this operation is involutive: swapping the roles of $a$ and $b$, and applying it to the line $S(S(\ell,O),b)$, takes one back to line $\ell$.

Following Section \ref{hyperbolic}, one can describe the hyperbolic case of
this construction in $\CP^2$;  the result is (a complex version of) the central projection in Figure \ref{projective}. 

\begin{figure}[hbtp]
\centering
\includegraphics[height=2.7in]{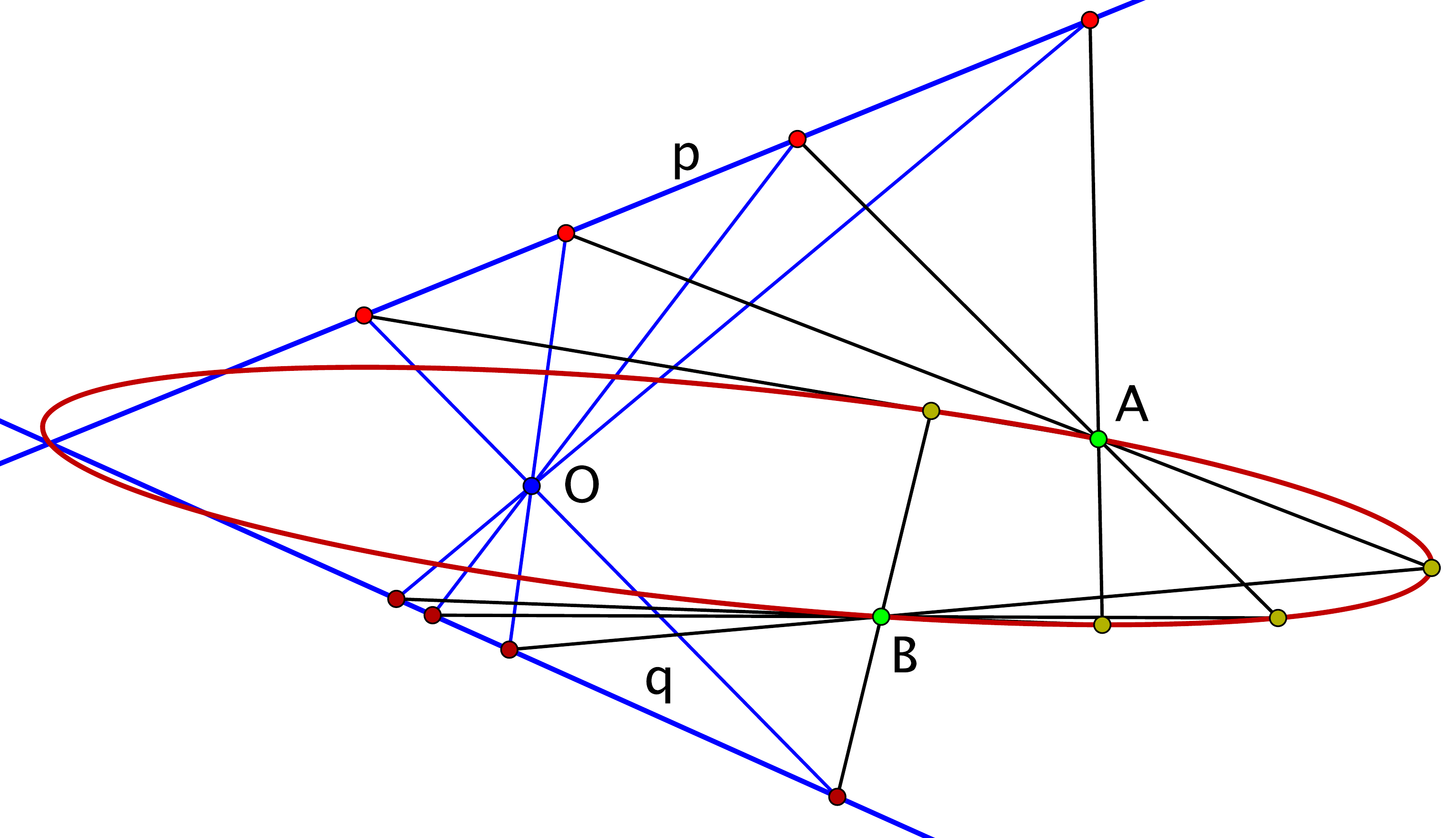}
\caption{The Braikenridge-Maclaurin construction of a conic.}
\label{conic}
\end{figure}

Recall the Braikenridge-Maclaurin construction of a conic depicted in Figure \ref{conic}; see \cite{Mi} for the history of this result.

Fix two lines, $p$ and $q$, and three points, $O, A$ and $B$. Identify $p$ with the pencil 
of lines through point $A$, and $q$  with the pencil of lines through point $B$.
The central projection $\varphi_O:p \to q$  induces a projective transformation between the two pencils of lines. Then the locus of intersection points of the corresponding lines from these pencils is a conic.

One can use the skewer version of the Braikenridge-Maclaurin construction to define a line analog of a conic. Start with five lines $O,p,q,A,B$. For each line $\ell \in {\cal N}_p$, we have the corresponding line $m=S(S(\ell,O),q) \in {\cal N}_q$. Then the 2-parameter family of lines
$$
S(S(\ell,A),S(m,B)),\ \ \ell \in {\cal N}_p
$$
is a skewer analog of a conic. In the hyperbolic case, this set is identified with a conic in $\CP^2$.

\section{Sylvester Problem} \label{Sylv}
Given a finite set $S$ of points in the plane, assume that the line through every pair of points in $S$ contains at least one other point of $S$. J.J.Sylvester asked in 1893 whether $S$ necessarily consists of collinear points. See \cite{BM} for the history of this problem and its generalizations.

In $\R^2$, the Sylvester question has an affirmative answer (the Sylvester-Galai theorem), but in $\C^2$ one has a counter-example: the 9 inflection points of a cubic curve (of which at most three can be real, according to a theorem of Klein), connected by 12 lines.

Note that the dual Sylvester-Galai theorem holds as well: if a finite collection of pairwise non-parallel lines in $\R^2$ has the property that through the intersection point of any two lines there passes at least one other line, then all the lines are concurrent.

The skewer version of the Sylvester Problem concerns a finite collection $S$ of pairwise skew lines in space such that the skewer of any pair intersects at least one other line at right angle. We say that $S$ has the skewer Sylvester property. 
The question is whether a collection of lines with the skewer Sylvester property necessarily consists of lines that share a skewer. 

\begin{theorem} \label{Sylvth}
The skewer version of the Sylvester-Galai theorem holds in the elliptic and Euclidean geometries, but fails in the hyperbolic one.
\end{theorem}

\proof In the elliptic case, we argue as in Section \ref{spherical}. A collection of lines becomes two collections of points,  in $\RP^2_-$ and in $\RP^2_+$.  The skewer Sylvester property implies that each of these  sets enjoys the property that the line through a pair of points contains another point, and one applies the Sylvester-Galai theorem to each of the two sets.

In the hyperbolic case, we argue as in Section \ref{hyperbolic}. Let $a_1,\ldots,a_9$ be the nine inflection points of a cubic curve in $\CP^2$, and let $b_1,\ldots,b_{12}$ be the respective lines (the counterexample to the complex Sylvester-Galai theorem). Let $b_1^*,\ldots,b_{12}^*$ be the polar dual points. As described in Section \ref{hyperbolic}, the points $a_i$ correspond to nine lines in $H^3$, and the points $b_j^*$ to their skewers. We obtain a collection of nine lines that has the skewer Sylvester property but does not possess a common skewer.

In the intermediate case of $\R^3$, the following argument is due to V. Timorin (private communication).

The approach is the same as in Section \ref{Eucpict}. It follows from the discussion there that if three lines in $\R^3$ share a skewer then their intersections with the plane at infinity $H$ are collinear. 

Let $L_1,\ldots, L_n$ be a collection of  lines enjoying the skewer Sylvester property. Then, by the Sylvester-Galai theorem in $H$, the points $q(L_1),\ldots, q(L_n)$  are collinear. This means that the lines $L_1,\ldots, L_n$ lie in parallel planes, say, the horizontal ones. 

Consider the vertical projection of these lines. We obtain a finite collection of non-parallel lines such that through the intersection point of any two there passes at least one other line. By the dual Sylvester-Galai theorem, all these lines are concurrent. Therefore the horizontal lines in $R^3$ share a vertical skewer. 
\proofend

\section{Pappus revisited} \label{Papprev}

In this section we prove Theorem \ref{othPapp}. This computational proof is joint with R. Schwartz.

As before, it suffices to establish the hyperbolic version of Theorem \ref{othPapp}.
We use the approach to 3-dimensional hyperbolic geometry, in the upper half-space model, developed by Fenchel \cite{Fe}; see also \cite{Iv,Ma}. The relevant features of this theory are as follows.

To a line $\ell$ in $H^3$, one assigns the reflection in this line, an orientation preserving isometry of the hyperbolic space,  an element of the group $PGL(2,\C)$. One can lift it to a matrix ${M_{\ell}} \in GL(2,\C)$, defined up to a complex scalar. Since reflection is an involution, one has ${\Tr} (M_{\ell})=0$. More generally, a traceless matrix $M \in GL(2,\C)$ is called a {\it line matrix}; it satisfies $M^2 = -\det(M) E$ where $E$ is the identity matrix.

The skewer relations translate to the language of matrices as follows:
\begin{itemize}
\item two lines $\ell$ and $n$ intersect at  right angle if and only if ${\Tr} (M_{\ell} M_n)=0$;
\item the skewer of two lines $\ell$ and $n$ corresponds to the commutator $[M_{\ell},M_n]$;
\item three lines $\ell,m,n$ share a skewer if and only if the matrices $M_{\ell}, M_m$, and $M_n$ are linearly dependent.
\end{itemize}

Likewise, one assigns matrices to points. The reflection in a point $P$ is an orientation-reversing isometry of $H^3$; one assigns to it a matrix $N_P$ in $GL(2,\C)$,  defined up to a real scalar, with $\det N_P >0$ and satisfying $N_P {\overline N_P} = -\det(N_P) E$, where bar means the entry-wise complex conjugation of a matrix. Such matrices are called {\it point matrices}. 

Equivalently, point matrices $N$ satisfy $n_{22}=-{\bar n_{11}},  n_{12}\in\R, n_{21}\in\R$, that is, the real part of $N$ is a traceless matrix, and the imaginary part is a scalar matrix. It is convenient to normalize so that the imaginary part is $E$, and then $N$ can be though of as a real 3-vector consisting of three entries of the real part of $N$.

Incidence properties translate as follows:
\begin{itemize}
\item a point $P$ lies on a line $\ell$ if and only of $M_{\ell} N_P = N_P {\overline M_{\ell}}$;
\item three points are collinear if and only if the respective point matrices are linearly dependent (equivalently, over $\R$ or $\C$).
\end{itemize}

We need a formula for a line matrix corresponding to the line through two given points. Let $N_1$ and $N_2$ be point matrices corresponding to the given points. Then the desired line matrix $M\in GL(2,\C)$ satisfies the system of linear equations
\begin{equation} \label{linepoint}
M N_1 = N_1 \overline {M},\ M N_2 = N_2 {\overline M},\ {\Tr} (M) =0.
\end{equation}
This system is easily solved and it defines $M$ up to a factor (we do not reproduce the explicit formulas here). 

With these preliminaries, the proof proceeds in the following steps. 
\begin{enumerate}
\item Start with two triples of linearly dependent point matrices, corresponding to the triples of points $A_1,A_2,A_3$ and $B_1,B_2,B_3$.
\item Compute the line matrices, corresponding to the lines $(A_1 B_2)$ and $(A_2 B_1)$, $(A_2 B_3)$ and $(A_3 B_2)$, and $(A_3 B_1)$ and $(A_1 B_3)$ by solving the respective systems (\ref{linepoint}).
\item Compute the commutators of these three pairs of line matrices.
\item Check that the obtained three matrices are linearly dependent.
\end{enumerate}

We did these computations in Mathematica. Since a line matrix is traceless, it can be viewed as a complex 3-vector, and the last step consists in computing the determinant made by three 3-vectors. The result of this last computation was zero (for arbitrary initial point matrices) which proves the theorem.

\begin{remark}
{\rm Theorem \ref{othPapp} can be restated somewhat similarly to  Theorem \ref{skPappus}. Given two skew lines $L$ and $M$, consider the 1-parameter family  of lines ${\cal F}(L,M)$ consisting of the lines that pass through a point $A\in L$ and orthogonal to the plane spanned by point $A$ and line $M$. Likewise, one has the 1-parameter family of lines ${\cal F}(M,L)$.  
These families, ${\cal F}(L,M)$ and ${\cal F}(M,L)$, replace the 2-parameter families of lines ${\cal N}_L$ and ${\cal N}_M$ in the formulation of Theorem \ref{skPappus}, and yield Theorem \ref{othPapp}. 
}
\end{remark}

\begin{remark}
{\rm F. Bachmann \cite{Ba} developed an approach to 2-dimensional geometry (elliptic, Euclidean, and hyperbolic) based on  the notion of reflection and somewhat similar to Fenchel's approach to 3-dimensional hyperbolic geometry \cite{Fe}. Namely, to a point $P$ there corresponds the reflection $\sigma_P$ in this point, and to a line $\ell$ -- the reflection $\sigma_{\ell}$ in this line. The incidence relation $P \in \ell$ is expressed as $\sigma_P \sigma_{\ell} = \sigma_{\ell} \sigma_P$. Two lines, $\ell$ and $m$, are orthogonal if and only if  $\sigma_{\ell} \sigma_m = \sigma_m \sigma_{\ell}$. More generally, one has a system of axioms of plane geometry in terms of involutions in the group of motions. At the present writing, it is not clear how to deduce the Correspondence principle using 
this approach. 
}
\end{remark}


\begin{thebibliography}{99}

\bibitem{AP1} N. A'Campo, A. Papadopoulos. {\it Transitional geometry.} in {\it Sophus Lie and Felix Klein: The Erlangen Program and its Impact in Mathematics and Physics}. L. Ji and A. Papadopoulos. ed.,  
IRMA Lect.  Math. Theor. Physics, 23. European Math. Soc., Z\"urich, 2015.


\bibitem{Ai} F. Aicardi. {\it Projective geometry from Poisson algebras.} J. Geom. Phys. {\bf 61} (2011), 1574--1586.


\bibitem{Ar} V. Arnold. {\it Lobachevsky triangle altitude theorem as the Jacobi identity in the Lie algebra of quadratic forms on symplectic plane}. J. Geom. Phys. {\bf 53} (2005), 421--427.

\bibitem{Ba} F. Bachmann. {\it Aufbau der Geometrie aus dem Spiegelungsbegriff.} Springer-Verlag, Berlin-G\"ottingen-Heidelberg, 1959.

\bibitem{BM} P. Borwein, W. O. J. Moser. {\it A survey of Sylvester's problem and its generalizations.} Aequationes Math. {\bf 40} (1990), 111--135. 

\bibitem{CS} S.-J. Chang, K. Shi. {\it Billiard systems on quadric surfaces and the Poncelet theorem.}  J. Math. Phys. {\bf 30} (1989),  798--804.

\bibitem{CR1} J. Conway, A. Ryba. {\it The Pascal mysticum demystified.} Math. Intelligencer {\bf 34} (2012), no. 3, 4--8. 

\bibitem{CR2} J. Conway, A. Ryba. {\it Extending the Pascal mysticum.} Math. Intelligencer {\bf 35} (2013), no. 2, 44--51.

\bibitem{Coo} J. L.  Coolidge. {\it A treatise on the circle and the sphere.} Reprint of the 1916 edition. Chelsea Publ. Co., Bronx, N.Y., 1971.

\bibitem{Co} H. S. M. Coxeter. {\it The inversive plane and hyperbolic space.}
 Abh. Math. Sem. Univ. Hamburg {\bf 29} (1966), 217--242.
 
\bibitem{Do} I. Dolgachev. {\it Classical algebraic geometry. A modern view.} Cambridge Univ. Press, Cambridge, 2012.
 
 \bibitem{DR} V.  Dragovi\'c, M.  Radnovi\'c. {\it Poncelet porisms and beyond. 
Integrable billiards, hyperelliptic Jacobians and pencils of quadrics.}  Birkh\"auser/Springer, Basel, 2011.

\bibitem {Fe} W. Fenchel. {\it Elementary geometry in hyperbolic space.} Walter de Gruyter, Berlin, 1989.

\bibitem {Fl} L. Flatto. {\it Poncelet's theorem.} Amer. Math. Soc., Providence, RI, 2009.

\bibitem {GW} H. Gluck, F. Warner. {\it Great circle fibrations of the three-sphere.} Duke Math. J. {\bf 50} (1983), 107--132.

\bibitem {Gr} B. Gr\"unbaum. {\it Configurations of points and lines}.  Amer. Math. Soc., Providence, RI, 2009.

\bibitem{HC} D. Hilbert, S.  Cohn-Vossen. {\it Geometry and the imagination.} Chelsea Publishing Company, New York, N. Y., 1952.

\bibitem {Ho} P. Hooper. {\it From Pappus' theorem to the twisted cubic.} Geom. Dedicata {\bf 110} (2005), 103--134. 

\bibitem{Iv} N. Ivanov. {\it Arnol'd, the Jacobi identity, and orthocenters.}
Amer. Math. Monthly {\bf 118} (2011),  41--65.

\bibitem{Je} C. M. Jessop. {A treatise on the line complex}. Cambridge Univ. Press, Cambridge, 1903.

\bibitem{KS} B. Konopelchenko, W. Schief. {\it Menelaus' theorem, Clifford configurations and inversive geometry of the Schwarzian KP hierarchy}. J. Phys. {\bf A 35} (2002), 6125--6144. 

\bibitem {Ma} A. Marden. {\it Outer circles. An introduction to hyperbolic 3-manifolds.} Cambridge Univ. Press, Cambridge, 2007. 

\bibitem {Mi} S. Mills. {\it Note on the Braikenridge-Maclaurin theorem.} Notes and Records Roy. Soc. London {\bf 38} (1984), 235--240.

\bibitem{Mo} F. Morley. {\it On a regular rectangular configuration of ten lines}. Proc. London Math. Soc. s1-29 (1897), 670--673.

\bibitem{Mo1} F. Morley. {\it The Celestial Sphere.} Amer. J. Math. {\bf 54} (1932),  276--278.

\bibitem{Mo2} F. Morley, F. V. Morley. {\it Inversive geometry}. G. Bell \& Sons, London, 1933.


\bibitem{Pa} {\it Strasbourg master class on geometry.} 
A. Papadopoulos ed. IRMA Lect.  Math. Theor. Physics, 18. European Math. Soc., Z\"urich, 2012.

\bibitem {PW} H. Pottmann, J. Wallner. {\it Computational line geometry.} Springer-Verlag, Berlin, 2001.

\bibitem {Ri} H. Richmond.  {\it A chain of theorems for lines in space.} J. London Math. Soc. {\bf 16} (1941), 108--112.

\bibitem{RG} J. Richter-Gebert. {\it Perspectives on projective geometry. A guided tour through real and complex geometry.} Springer, Heidelberg, 2011.

\bibitem{SK} W. Schief, B. Konopelchenko. {\it A novel generalization of Clifford's classical point-circle configuration. Geometric interpretation of the quaternionic discrete Schwarzian Kadomtsev-Petviashvili equation.} Proc. R. Soc. Lond. Ser. A Math. Phys. Eng. Sci. {\bf 465} (2009), 1291--1308.

\bibitem {Sc1} R. Schwartz. {\it The pentagram map}. Experiment. Math. {\bf 1} (1992), 71--81.

\bibitem {Sc2} R. Schwartz. {\it Pappus' theorem and the modular group}. Inst. Hautes \'Etudes Sci. Publ. Math. No. 78 (1993), 187--206 (1994).

\bibitem {ST} R. Schwartz, S. Tabachnikov. {\it Elementary surprises in projective geometry}. Math. Intelligencer {\bf 32} (2010), no. 3, 31--34. 

\bibitem {Sk} M. Skopenkov. {\it Theorem about the altitudes of a triangle and the Jacobi identity} (in Russian). Matem. Prosv., Ser. 3 {\bf 11} (2007), 79--89.

\bibitem{St} E. Study. {\it Geometrie der Dynamen}. B. G. Teubner, Leipzig, 1903 

\bibitem{To} T. Tomihisa. {\it Geometry of projective plane and Poisson structure.} J. Geom. Phys. {\bf 59} (2009),   673--684.

\bibitem{Ve} A. Veselov. {\it Confocal surfaces and integrable billiards on the sphere and in the Lobachevsky space.} J. Geom. Phys. {\bf 7} (1990),  81--107. 

\end{thebibliography}
\end{document}